\documentclass[12pt]{amsart}
 
 \usepackage{graphicx,hyperref}

\numberwithin{equation}{section}
 \usepackage{bbold}
 \usepackage{amsmath}
\usepackage{amsthm}
\usepackage{amscd}
\usepackage{amssymb}
\usepackage{comment}

\theoremstyle{plain}
\newtheorem*{theorem*}{Theorem}
\newtheorem{theorem}[equation]{Theorem}

\newtheorem{lemma}[equation]{Lemma}
\newtheorem{corollary}[equation]{Corollary}

 \theoremstyle{definition}
\newtheorem{definition}[equation]{Definition}
\newtheorem{remark}[equation]{Remark}

\newcommand{\G}{{\mathcal G}}
 \renewcommand\H{{\mathcal H}}

\newcommand{\R}{\mathbf{R}}

\newcommand{\al}{\alpha}
\newcommand{\ant}{\operatorname{Ant}}
\newcommand{\area}{\operatorname{Area}}

\newcommand{\be}{\beta}
\newcommand{\cat}{\operatorname{CAT}}
\newcommand{\cone}{C}
\newcommand{\ctits}{C_T}
\newcommand{\de}{\delta}
\newcommand{\defeq}{:=}
\newcommand{\diam}{\operatorname{diam}}
\newcommand{\eps}{\epsilon}

\newcommand{\ga}{\gamma}
\newcommand{\Ga}{\Gamma}

\newcommand{\id}{\operatorname{id}}
\newcommand{\im}{\operatorname{Im}}

\newcommand{\lf}{\operatorname{lf}}
\newcommand{\lra}{\longrightarrow}

\newcommand{\ol}{\overline}

\newcommand{\ra}{\rightarrow}

\newcommand{\si}{\sigma}
\newcommand{\Si}{\Sigma}
\newcommand{\supp}{\operatorname{supp}}
\renewcommand{\th}{\theta}
\newcommand{\tits}{\partial_T}

\begin{document}

\title{Quasiflats in $\cat(0)$ complexes}
\author{Mladen Bestvina, Bruce Kleiner, and Michah Sageev}
\thanks{This research was supported by NSF grants DMS-1308178 and DMS-1405899}
\maketitle 

\begin{abstract}
We show that if $X$ is a piecewise Euclidean $2$-complex with
a cocompact isometry group, then every $2$-quasiflat in $X$
is at finite Hausdorff distance from a subset $Q$ which is 
locally flat outside a compact set, and asymptotically conical.
\end{abstract}

\section{Introduction}

In a number of rigidity theorems for quasi-isometries, an
important step is to determine the structure of individual
quasi-flats; this is then  used to restrict the behavior of 
quasi-isometries, often by exploiting the pattern of asymptotic
incidence of the quasiflats 
\cite{mostow,kaple,klle,eskinfarb,eskin,behrstockkleinerminskymosher}.
In this paper, we study $2$-quasiflats in $\cat(0)$ $2$-complexes,
and show that they have a very simple asymptotic structure:

\begin{theorem}
\label{main}
Let $X$ be a proper, piecewise Euclidean, $\cat(0)$ $2$-complex
with a cocompact isometry group.  Then every $2$-quasiflat 
$Q\subset X$ lies at finite Hausdorff distance from
a subset $Q'\subset X$ which is locally flat (i.e. locally
isometric to $\R^2$) outside a compact set.  
\end{theorem}  

This result, and more refined statements appearing in later sections,
are applied to $2$-dimensional right-angled Artin groups in \cite{raag}.  
The main application is to show that if $X,\,X'$ are
the standard $\cat(0)$ complexes of $2$-dimensional right-angled
Artin groups, then any quasi-isometry $X\ra X'$
between them must map flats to within finite Hausdorff distance of flats.

The strategy for proving Theorem \ref{main} is to replace the  quasiflat
$Q$ with a canonical object that has more rigid structure. 
To that end, we first associate an element 
$[Q]$ of the locally finite homology group $H_2^{\lf}(X)$, and
then show that the support set $\supp([Q])$ of $[Q]$ -- the set
of points $x\in X$ such that the 
induced homomorphism $H_2^{\lf}(X)\ra H_2(X,X\setminus \{x\})$
is nontrivial on $[Q]$ -- is at bounded Hausdorff distance from $Q$.
The support set $Q'\defeq\supp([Q])$ behaves much like a minimizing locally finite
cycle, and this leads to asymptotically rigid behavior, in 
particular asymptotic flatness.

\bigskip
{\bf Remarks.}
\begin{enumerate}
\item Support sets were used implicitly in \cite{klle},
and also in \cite{xie}.  

\item The paper \cite{klla}, which may
be viewed as a more sophisticated version of the
results presented here,  exploits similar geometric 
ideas in asymptotic cones, to study $k$-quasiflats in 
$\cat(0)$ spaces which have no $(k+1)$-quasiflats.

\item Many of the results of this paper (though not 
Theorem \ref{main} itself) can be adapted 
 to $n$-quasiflats in $n$-dimensional $\cat(0)$ complexes.

\item One may use the results in this paper to give a new
proof that quasi-isometries between Euclidean buildings map flats
to within uniform Hausdorff distance of flats \cite{klle}.
This then leads to a (partly) 
different proof of rigidity of quasi-isometries
between Euclidean buildings.
\end{enumerate}

\tableofcontents

\section{Preliminaries}

\subsection{$\cat(\kappa)$ spaces}
We recall some standard facts, and fix notation.
We refer the reader to \cite{bridsonhaefliger,klle} for more detail.
Our notation and conventions are consistent with \cite{klle}.

Let $X$ be a $\cat(0)$ space.

If $x,y\in X$, then $\ol{xy}\subset X$
denotes the geodesic segment with endpoints $x,y$. If $p\in X$,
we let $\angle_p(x,y)$ denote the angle between $x$ and $y$ at $p$.
This
induces a pseudo-distance on $X\setminus\{p\}$.  By collapsing subsets
of zero diameter and completing, we obtain   
the space of directions $\Si_pX$, which is a $\cat(1)$ space.
The quotient map yields the logarithm $\log_p:X\setminus\{p\}\ra \Si_pX$;
it associates to $x\in X\setminus \{p\}$ the direction at $p$ of the 
geodesic segment $\ol{px}$.  
The tangent cone at $p$, denoted
$\cone_pX$,  is a $\cat(0)$ space isometric to the cone over $\Si_pX$.

Given two constant (not necessarily unit) speed rays
$\ga_1,\ga_2:[0,\infty)\ra X$,  their distance is defined to be
$$
\lim_{t\ra\infty}\;\frac{d(\ga_1(t),\ga_2(t))}{t}.
$$
This defines a pseudo-distance on the set of constant speed rays 
in $X$; the metric space obtained by collapsing zero diameter
subsets is the {\em Tits cone} of $X$, denoted
$\ctits X$.  The Tits cone is isometric to
the Euclidean cone over the
Tits boundary $\tits X$. 
For every $p\in X$, there are natural logarithm maps:
$$
\log_p:X\ra \cone_pX,\quad \log_p:\ctits X\ra X, 
$$
$$
\log_p:X\setminus\{p\}\ra \Si_pX,\quad \log_p:\tits X\ra \Si_pX.
$$

\bigskip
\bigskip
\begin{definition}
\label{defantipodal}
If $Z$ is a $\cat(1)$ space, $Y\subset Z$, and $z\in Z$, then 
the {\em antipodal set of $z$ in $Y$},  is
$$
\ant(z,Y)\defeq \{y\in Y\mid d(z,y)=\pi\}.
$$
Recall that by our definition, every $\cat(1)$ space has diameter $\leq \pi$.
\end{definition}

If $X$ is a $\cat(0)$ complex and $p,x\in X$ are distinct points, $Y\subset\Si_xX$,
 then the antipodal set $\ant(\log_xp,Y)$ is the set of directions in $Y$
which are tangent to extensions of  the geodesic segment $\ol{px}$
beyond $x$.

\bigskip
\subsection{Locally finite homology}

Let $Z$ be a topological space. We recall that  the $k^{th}$
locally finite (singular) chain group $C^{\lf}_k(Z)$ is the 
collection of (possibly infinite) formal sums of singular $k$-simplices,
such that for every compact subset $Y\subset Z$, only finitely
many nonzero terms are contributed by singular simplices whose
image intersects $Y$.  The usual boundary operator yields a 
well-defined chain complex $C_*^{\lf}(Z)$; its homology is
the {\em locally finite homology of $Z$}.

Suppose  $K$ is a simplicial complex. Then there is a simplicial
version of the locally finite chain complex -- the locally finite
simplicial chain complex -- defined by taking (possibly infinite)
formal linear combinations of oriented simplices of $K$, where 
every simplex $\si$ of $K$ touches only finitely many simplices with
nonzero coefficients.
The usual proof that simplicial homology is isomorphic to
singular homology gives an isomorphism between the locally finite
simplicial homology of $K$, and the locally finite homology of
its geometric realization $|K|$, when $K$ is locally finite
\cite[3.H, Exercise 6]{hatcher}. 

The {\em support set} of $\si\in H_k^{\lf}(Z)$ is the subset
$\supp(\si)\subset Z$ consisting of the points $z\in Z$
for which the inclusion homomorphism
$$
H_k^{\lf}(Z)\ra H_k(Z,Z\setminus\{z\})
$$
is nonzero on $\si$.  This is a closed subset when $Z$ is Hausdorff.  

Now suppose $K$ is an $n$-dimensional locally finite simplicial
complex, with polyhedron $Z$.  Then the simplicial chain groups 
$C_k^{\lf}(K)$ vanish for $k>n$, and hence $H_n^{\lf}(Z)$
is isomorphic to the group of locally finite simpicial
$n$-cycles $Z_n^{\lf}(K)$.  The support set of a locally finite simplicial
$n$-cycle $\si\in Z_n^{\lf}(Z)$ is the union of the closed $n$-simplices
having  nonzero coefficient in $\si$, as follows from excision.

\section{Locally finite homology and support sets}
The key results in this section  are the geodesic extension property
of Lemma \ref{lemextend}, and the asymptotic conicality result for 
support sets with quadratic area growth, in Theorem 
\ref{thmquadraticgrowthsupportsetstructure}.  We remark that most of the statements
(and proofs) in this section extend with minor modifications
to supports of $n$-dimensional locally finite homology classes in  $n$-dimensional 
$\cat(0)$ complexes.

In this section $X$ will be a proper, piecewise
Euclidean, $\cat(0)$ $2$-complex.

\subsection{The geodesic extension property and metric monotonicity}
The fundamental property of support sets is the extendability of geodesics:

\begin{lemma}
\label{lemextend}
Suppose $\si\in H_2^{\lf}(X)$, and let $S\defeq \supp(\si)\subset X$ be the 
support of $\si$.
If $p\in X$, and $x\in S$,  the geodesic segment $\ol{px}$ 
may be prolonged to a ray in $S$:
there is a ray $\ol{x\xi}\subset S$ with fits together with
$\ol{px}$ to form a ray $\ol{p\xi}$.
\end{lemma}
\proof
Let $\ga:[0,L]\ra X$ be the unit speed parametrization of $\ol{px}$,
and let  $\hat\ga:I\ra X$ be a maximal extension of  $\ga$ such that
$\hat\ga(I\setminus [0,L])\subset S$, where $I$ is an interval
contained in $[0,\infty)$.  Since $S$ is  a closed subset of the complete
space $X$, either $I=[0,R]$ for some $R<\infty$, or $I=[0,\infty)$.

Suppose $I=[0,R]$ for $R<\infty$, and let $y\defeq \hat\ga(R)$.
Consider the closed ball $B\defeq\ol{B(y,r)}$, where $r$ is small enough that 
$B$ is isometric
to the $r$-ball in the tangent cone $C_yX$.  Note that this implies that $S\cap B$ is also
a cone.  Let $\si=[\si_B+\tau]$, where
$\si_B\in C_2^{\lf}(X)$ is  carried by  $B$ (and is therefore a finite 
$2$-chain),  $\tau\in C_2^{\lf}(X)$ is 
carried by  $X\setminus B(y,r)$, and $\partial \si_B=-\partial \tau$ is 
carried by   $\partial B\cap S$.
Consider the singular chain $\mu$ obtained by coning off $\partial\si_B$
at $p$.  Then $\partial \mu=\partial \si_B$, so the contractibility of $X$
implies that $\mu$ is homologous to $\si_B$ relative to $\partial \mu$.
Thus $\mu+\tau$ belongs to the homology class of  $\si$.  Therefore $y$ lies
in the carrier of $\mu$, for otherwise $\mu+\tau$ would be carried by 
$X\setminus \{y\}$, contradicting the fact that $y\in\supp(\si)$. 
Thus there is a point $z\in \partial B\cap S$
such that the segment $\ol{pz}$ passes through $y$.  Since $B\cap S$ is a
cone, we have $\ol{yz}\subset S$.   This implies that $\hat\ga$ is not a maximal
extension, which is a contradiction.

Another way to argue the last part of the proof is to observe that $\si_B$
projects under $\log_y:X\setminus\{y\}\ra\Si_yX$
 to a nontrivial $1$-cycle $\eta$ in $\Si_yX$.  Therefore,
there must be a direction $v\in \Si_yS$ making an angle $\pi$ with
$\log_yp$, since otherwise $\eta$ would lie in the open ball
of radius $\pi$ centered at $\log_yp$, which is contractible.  Then
$\hat\ga$ may be extended in the direction $v$, which contradicts
the maximality of $\hat\ga$. 
\qed

\begin{remark}
The geodesic extension property has a flavor similar to convexity, but note that
support sets need not be convex.  To obtain an example, let $Z$ be
the union of two disjoint circles $Y_1,\,Y_2$ of length
$2\pi$ with a geodesic segment of length $<\pi$ (so $Z$ is a
``pair of glasses''), and let $X$ be the Euclidean cone over $Z$.  Then 
cone over $Y_1\cup Y_2$ is a support set, but is not convex.
\end{remark}

\begin{corollary}[Monotonicity and lower density bound]
\label{cormonotonicitydensity}
Suppose $\si\in H_2^{\lf}(X)$ and $S\defeq \supp(\si)$.

1. (Metric monotonicity)  
For all $0<r\leq R$, $p\in X$, if $\Phi:X\ra X$ is the map
which contracts points toward $p$ by the factor $\frac{r}{R}$, then
\begin{equation}
\label{eqnmetricmonotonicity}
B(p,r)\cap S\;\subset\;\Phi(B(p,R)\cap S).
\end{equation}

2. (Monotonicity of density)  For all $0\leq r\leq R$, 
\begin{equation}
\label{eqndensitymonotone}
\frac{\area(B(p,r)\cap S)}{r^2}\leq \frac{\area(B(p,R)\cap S)}{R^2}.
\end{equation}

3. (Lower density bound) For all $p\in S$, $r>0$, 
\begin{equation}
\label{eqnlowerdensity}
\area(B(p,r)\cap S)\geq  \pi\, r^2,
\end{equation}
with equality only if $B(p,r)\cap S$ is isometric to an $r$-ball in $\R^2$.

Here $\area(Y)$ refers to $2$-dimensional Hausdorff measure, 
which is the same as Lebesgue measure (computed by summing over the
intersections with $2$-dimensional faces).
\end{corollary}

\begin{remark}
Since the map $\Phi$ in assertion 1 has Lipschitz constant
$\frac{r}{R}$, the inclusion
(\ref{eqnmetricmonotonicity})  can be viewed as a much stronger version of the 
usual monotonicity formula for minimal submanifolds in nonpositively curved spaces, which
corresponds to (\ref{eqndensitymonotone}).
\end{remark}

{\em Proof of Corollary \ref{cormonotonicitydensity}.}
(\ref{eqnmetricmonotonicity}) follows from Lemma \ref{lemextend}.

Assertion 2 follows from assertion 1 and the fact that $\Phi$ has
Lipschitz constant $\frac{r}{R}$.

If $p\in S$, then $\si$ determines a nonzero class $\Si_p\si\in H_1(\Si_pX)$,
by the composition
$$
H_2(X,X\setminus\{p\})\stackrel{\partial}{\lra}H_1(X\setminus\{p\})\stackrel{\log_{\Si_pX}}{\lra}
H_1(\Si_pX).
$$
Since $\Si_pX$ is a $\cat(1)$ graph, $\supp(\Si_p\si)$ contains a cycle of length 
at least $2\pi$.  If $r>0$ is small, then $B(p,r)\cap S$ is isometric to a
cone of radius $r$ over $\supp(\Si_p\si)$, and therefore has area at least
$\pi\,r^2$.  Now (\ref{eqndensitymonotone}) implies (\ref{eqnlowerdensity}).  
Equality in (\ref{eqnlowerdensity})
implies that $\supp(\Si_p\si)$ is a circle of length $2\pi$, $B(p,r_0)\cap S$ is isometric
to an $r_0$-ball in $\R^2$ for small $r_0>0$, and that the contraction map $\Phi$
is similarity.   This implies 3.
\qed

\bigskip
The corollary implies that the ratio 
$$
\frac{\area(B(p,r)\cap S)}{r^2}
$$
has a (possibly infinite) limit $\bar A$
as $r\ra\infty$, which   is clearly independent of the basepoint.
When it is finite we say that $\si$ has {\em quadratic growth}.  In this case,
Corollary \ref{cormonotonicitydensity} implies that 
\begin{equation}
\label{eqnquadraticgrowth}
\frac{\area(B(p,r)\cap S)}{r^2}\leq \bar A
\end{equation}
 for all $p\in X$, $r>0$.

\bigskip
\bigskip

\subsection{Asymptotic conicality}
We will use Lemma \ref{lemextend} and Corollary \ref{cormonotonicitydensity}
to see that quadratic growth support sets are asymptotically conical, provided
the $\cat(0)$ $2$-complex $X$ satisfies a mild additional condition.  To see
why an additional assumption is needed, consider a piecewise
Euclidean $\cat(0)$ $2$-complex $X$
homeomorphic to $\R^2$, whose singular set consists of a
sequence of cone points  $\{p_i\}$  tending to infinity, 
where $\Si_{p_i}X$ is a circle
of length $2\pi+\th_i$, and $\sum_i\;\th_i<\infty$.  Then $X$ is the support
set of the locally finite fundamental class $[X]$ of the $2$-manifold $X$,
but it is not locally flat outside any compact subset of $X$.   

To exclude this kind of behavior, one would like to know, for instance,
that the cone angle  $2\pi$ is isolated among the set of cone angles of
points in $X$.  When dealing with general $\cat(0)$ $2$-complexes, one needs 
to know that if $p\in X$ and  $v\in \Si_pX$ is a 
direction whose antipodal set $\ant(v,\supp(\tau))$
in a $1$-cycle $\tau\in Z_1(\Si_pX)$ has small diameter, then $v$ is 
close to a suspension point of $\tau$.  This condition will hold automatically
if $X$ admits a cocompact group of isometries.  The precise 
condition we need is:

\begin{definition}
A family ${\mathcal F}$ of $\cat(1)$ graphs has {\em isolated suspensions} if for
every $\al>0$ there is a $\be>0$ such that if $\Ga\in {\mathcal F}$, 
$\tau\in Z_1(\Ga)$ is a $1$-cycle, $v\in \Ga$, and 
$$
\diam(\ant(v,\supp(\tau))<\be,
$$
then $\supp(\tau)$ is a metric suspension and $v$ lies at distance 
$<\al$ from a pole (i.e. suspension point) of $\supp(\tau))$.   
A $\cat(0)$ $2$-complex $X$ has {\em isolated suspensions} if the collection of
spaces of directions $\{\Si_xX\}_{x\in X}$ has isolated suspensions.
\end{definition}

\begin{remark}
\label{remfiniteisolated}
It follows from a compactness argument that any finite collection of $\cat(1)$
graphs has the isolated suspensions property.  In particular, any $\cat(0)$ $2$-complex
with a cocompact isometry group has the isolated suspension property.
\end{remark}

For the remainder of this section $X$ will be a piecewise Euclidean, proper
$\cat(0)$ $2$-complex with isolated suspensions.

\begin{theorem}
\label{thmquadraticgrowthsupportsetstructure}
Suppose $\si\in H_2^{\lf}(X)$ has quadratic area growth,
and $S\defeq\supp(\si)$.  Then for all $p\in X$ there is an $r_0<\infty$ such that

1. If 
$x\in S\setminus B(p,r_0)$, then  $S$ is locally isometric to a product
of the form $\R\times W$ near $x$, where $W$ is an $i$-pod (i.e. a cone
over a finite set).  In particular
$S$ is locally convex near $x$.

2.  The map $S\setminus B(p,r_0)\ra [r_0,\infty)$
given by the distance function from $p$ is a fibration with fiber homeomorphic
to a finite graph with all vertices of valence $\geq 2$.

3.  $S$ is asymptotically conical, in the following sense.
For every $p\in X$ and every $\eps>0$, there is an $r<\infty$ such that
if $x\in S\setminus B(p,r)$, then the angle (at $x$) between the geodesic
segment $\ol{xp}$ and the $\R$-factor of some local product splitting
of $S$ is $<\eps$.

4. If the area growth of $S$ is Euclidean, i.e.
$$
\frac{\area(B(p,r)\cap S)}{\pi r^2}\ra 1
$$
as $r\ra \infty$, then $S$ is a $2$-flat.
\end{theorem}

Before entering into the proof of this theorem,
we point out that the proof is driven by the following observation.
The locally finite cycle $\si$ is an area minimizing object in the 
strongest possible sense: any compact piece $\tau$ solves the Plateau
problem with boundary condition $\partial \tau$
(i.e. filling $\partial \tau$ with a least area chain); in fact, because of the 
dimension assumption, there is only one way to fill $\partial \tau$
with a chain.  Then we adapt
the standard monotonicity formula from minimal surface theory to see that 
the support set is asymptotically conical.  Roughly speaking the idea is that
the ratio 
$$
\frac{\area(B(p,r)\cap \supp(\si))}{r^2}
$$
is nondecreasing and bounded above, and hence has  limit as $r\ra \infty$. 
For large $r$, one concludes that the monotonicity inequality is nearly 
an equality, which leads to 2 of Theorem \ref{thmquadraticgrowthsupportsetstructure}.

{\em Proof of Theorem \ref{thmquadraticgrowthsupportsetstructure}.}
We begin with a packing estimate.

\medskip
\begin{lemma}
\label{lempackingbound}
For all $\eps>0$ there is an $N$ such that for all $r\geq 0$,
the intersection $B(p,r)\cap S$ does not
contain an $\eps r$-separated subset of cardinality greater than $N$.
\end{lemma}
\proof
Take $\eps<1$, and suppose the points $$x_1,\ldots,x_k\in B(p,r)\cap S$$
are $\eps r$-separated.  Then the collection
$$
\left\{\;B\left(x_i,\frac{\eps r}{2}\right)\cap S\right\}_{1\leq i\leq k}
$$
is disjoint, is contained in $B(p,2r)\cap S$, and by
assertion 2 of Corollary \ref{cormonotonicitydensity} it has area at least
$k\pi(\frac{\eps r}{2})^2$.  Thus (\ref{eqnquadraticgrowth}) implies
the lemma.
\qed

\bigskip
\begin{lemma}
\label{lemsmallant}
For all $\be>0$ there is an $r<\infty$ such
that if $x\in S\setminus B(p,r)$, then 
\begin{equation}
\label{eqndiameps}
\diam(\ant(\log_xp\;,\;\Si_xS))<\be.
\end{equation}\end{lemma}
\proof
The idea is that 
quadratic area growth
bounds the complexity of the support set from above, which implies
that on sufficiently large scales, it looks very much
 like a metric cone.   On the other hand, failure of (\ref{eqndiameps}) 
implies that there is a pair of rays leaving $p$
which coincide until $x$, and then
branch apart with an angle at least $\be$; when $x$ is far enough from
$p$, this will contradict the approximately
conical structure of $S$ at large scales.

Pick $\de,\mu>0$, to be determined later.  

By Lemma \ref{lempackingbound}
there is finite upper bound on the cardinality of an $\de r$-separated
subset sitting in $B(p,r)\cap S$, where $r$ ranges over $[1,\infty)$.
Let $N$ be the maximal such cardinality, which will be attained by some $\de r_0$-separated
subset $\{x_1,\ldots,x_N\}\subset B(p,r_0)\cap S$, for some $r_0$. 
Applying Lemma \ref{lemextend}, let $\ga_1,\ldots,\ga_N:[0,\infty)\ra X$
be constant speed geodesics emanating from $p$, such that 
$\ga_i(r_0)=x_i$, and $\ga_i(t)\in S$ for all $t\in[r_0,\infty)$,
$1\leq i\leq N$.   The  functions 
\begin{equation}
\label{eqnnondecrease}
t\mapsto \frac{d(\ga_i(t),\ga_j(t))}{t}
\end{equation}
are nondecreasing, and hence for all $r\in [r_0,\infty)$ the collection
$$
\ga_1(r),\ldots,\ga_N(r)
$$
is $\de r$-separated, and by maximality, it is therefore a $\de r$-net in $B(p,r)\cap S$
as well.  Using the monotonicity (\ref{eqnnondecrease}) again, we may find
$r_1\in [r_0,\infty)$ such that for all $1\leq i,j\leq N$, and every $r\in [r_1,\infty)$,

\begin{equation}
\label{eqnnearlimit}
\frac{d(\ga_i(r),\ga_j(r))}{r}+\mu>\lim_{t\ra\infty}\frac{d(\ga_i(t),\ga_j(t))}{t}\;.
\end{equation}

Now suppose $x\in S\setminus B(p,r_1)$, and $v_1,v_2\in\ant(\log_xp,\Si_xS)$
satisfy $\angle_x(v_1,v_2)\geq\be$.  The idea of the rest of the 
proof is  to invoke Lemma \ref{lemextend} to produce
two rays emanating from $p$ which agree until they reach
$x$, but then ``diverge at  angle at least $\beta$'';  since both rays
will be well-approximated by one of the $\ga_i$'s, their separation
behavior will contradict (\ref{eqnnearlimit}).

Let $r_2\defeq d(p,x)$.
By Lemma \ref{lemextend} we may prolong the segment $\ol{px}$
into two rays $\ol{p\xi_1},\ol{p\xi_2}$, such that $\log_{\Si_x}\xi_i=v_i$,
and $\ol{p\xi_i}\setminus B(p,r_2)\subset S$.  Let $\eta_1,\eta_2$
be the unit speed parametrizations of $\ol{p\xi_1}$ and $\ol{p\xi_2}$ respectively.
Applying triangle comparison, we may choose an $r_3\geq r_2$ such that
\begin{equation}
d(\eta_1(r_3),\eta_2(r_3))> r_3\cos\frac{\be}{2}.
\end{equation}
Pick $i,j$ such that 
$$
d(\ga_i(r_3),\eta_1(r_3))<\de r_3
\quad\mbox{and}\quad
d(\ga_j(r_3),\eta_2(r_3))<\de r_3.
$$
By triangle comparison, we have
$$
d(\ga_i(r_3),\ga_j(r_3))\geq d(\eta_1(r_3),\eta_2(r_3))-2\de r_3
>r_3\cos\frac{\be}{2}-2\de r_3
$$
while
$$
d(\ga_i(r_2),\ga_j(r_2))\leq d(\ga_i(r_2),\eta_1(r_2))+d(\eta_1(r_2),\eta_2(r_2))
+d(\eta_2(r_2),\ga_j(r_2))
$$
$$
\leq 2\de r_2,
$$
since $d(\eta_1(r_2),\eta_2(r_2))=0$.
On the other hand, by (\ref{eqnnearlimit})
$$
\mu>\frac{d(\ga_i(r_3),\ga_j(r_3))}{r_3}-\frac{d(\ga_i(r_2),\ga_j(r_2))}{r_2}
$$
$$
\geq \cos\frac{\be}{2}-4\de.
$$
When $\mu+4\de<\cos\frac{\be}{2}$ this gives a contradiction.
\qed

\bigskip
The Lemma together with the definition of isolated suspensions
implies parts 1 and 3 of Theorem \ref{thmquadraticgrowthsupportsetstructure}.
Part 4 follows from Lemma \ref{cormonotonicitydensity}.

To prove 2 of 
Theorem \ref{thmquadraticgrowthsupportsetstructure},
  we apply  the definition of isolated suspensions with $\al_0=\frac{\pi}{4}$ and let
$\be_0>0$ be the corresponding constant;   then we apply Lemma \ref{lemsmallant}
with $\be=\be_0$, and let $r_0$ be the resulting radius.  For each $x\in X\setminus B(p,r_0)$,
the space of directions $\Si_x S$ is a metric suspension,
and the direction $\log_xp\in\Si_xX$ makes an angle at most $\frac{\pi}{4}$ from a pole
of $\Si_xS$.

We call a point  $x\in S\setminus B(p,r_0)$ {\em singular} if its tangent cone is
not isometric to $\R^2$; thus singular points in $S\setminus B(p,r_0)$
have tangent cones of the form $\R\times W$, 
where $W$ is an $i$-pod with $i>2$, and the set of regular points forms an open
subset which carries the structure of a flat Riemannian manifold.  Using a partition of 
unity, we may construct a smooth  vector field $\xi$ on the regular
part of $S\setminus B(p,r_0)$ such that 

$\bullet$ $\xi(x)$ makes an angle at least $\frac{3\pi}{4}$
with $\log_xp$ at every regular point $x$.

$\bullet$ For each singular point $x\in S\setminus B(p,r_0)$ whose space of 
directions is the metric suspension of an
$i$-pod, if we decompose a small
neighborhood $B(x,\rho)\cap S$ into a union 
$$
C_1\cup\ldots\cup C_i,
$$
where the $C_j$'s are Euclidean half-disks of radius $\rho$ which 
intersect along a segment $\eta$ of 
length $2\rho$, then the restriction of $\xi$ to $C_j$ extends to a smooth
vector field $\xi_j$ on the manifold with boundary $C_j$, and $\xi_j(y)$
is a unit vector tangent to $\eta=\partial C_j$ for every $y\in \eta$.

\bigskip
Now a standard Morse theory argument using a reparametrization of the flow of $\xi$  implies
that 
$$
d_p:S\setminus B(p,r_0)\ra [r_0,\infty)
$$
is a fibration, and that the fiber is locally homeomorphic to an $i$-pod
near each point $x\in S\setminus B(p,r_0)$ whose space of directions is the
metric suspension of an $i$-pod.  Here $i\geq 2$.
\qed

\subsection{Asymptotic branch points}
The next result will be used when we consider support sets associated with
quasiflats.

\begin{lemma}
\label{lembranching}
Let $\si\in H_2^{\lf}(X)$ be a quadratic growth class with support $S$, 
pick $p\in X$, and let
$$
d_p:S\setminus B(p,r_0)\ra [r_0,\infty)
$$ be the fibration as in 2 of Theorem \ref{thmquadraticgrowthsupportsetstructure}.  
If the fiber has a branch point,  then for all $R<\infty$, the support set
$S$ contains an isometrically  embedded copy of an $R$-ball 
\begin{equation}
\label{eqnbr}
B_R\defeq B(z,R)\subset \R\times W,
\end{equation} 
where $W$ is an infinite
tripod, and $z\in \R\times W$ lies on the singular line.
\end{lemma}
\proof
Let $\pi:Y\ra S\setminus B(p,r_0)$ be the universal covering map.  Since
$S\setminus B(p,r_0)$ is homeomorphic to $\G\times [0,\infty)$, the covering
map $\pi$ is equivalent to the product of the universal covering $\tilde\G\ra \G$
with the identity map $[0,\infty)\ra[0,\infty)$. Since $\G$ contains a branch point,
we may find a proper embedding  $\phi:V\ra \tilde\G$ of a tripod $V$ into $\tilde \G$.
Consider the map $\psi$ given by the composition 
$$
V\times [0,\infty)\lra \tilde\G\times[0,\infty)\lra \G\times[0,\infty)\simeq S\setminus B(p,r_0).
$$
We may put a locally $\cat(0)$ metric on $V\times (0,\infty)$ by pulling back the 
metric from $S\setminus B(p,r_0)$.  For each of the three ``rays'' $\ga_i\subset V$
whose union is $V$, the metric on $\ga_i\times (0,\infty)$ is locally isometric
to a flat metric with geodesic boundary.  It follows from a standard argument
that if $y\in V\times (0,\infty)$ lies on the singular locus and $\psi(y)$ 
lies outside $B(p,r_0+R)$, then the $R$-ball in $V\times (0,\infty)$
is isometric to  $B_R$  as in (\ref{eqnbr}).  Since $\psi$ is a locally
isometric map of a $\cat(0)$ space into a $\cat(0)$ space, it is an isometric embedding.
\qed

\section{Quasi-flats in $2$-complexes}

In this section, $X$ will denote
 a piecewise flat, proper $\cat(0)$ $2$-complex with isolated 
suspensions.
\begin{theorem}
\label{thmqflataflat}
Let $Q\subset X$ be an $(L,A)$-quasiflat.  Then there is a 
nontrivial quadratic 
growth, locally finite homology class $\si\in H_2^{\lf}(X)$ whose
support set
$S\subset X$ is at Hausdorff distance at most $D=D(L,A)$ from $Q$, with
the following property.

1.  For every $p\in X$, there is an $r_0\in [0,\infty)$ such that
 $S\setminus \ol{B(p,r_0)}$ is  locally isometric to $\R^2$.

2.   $S$ is asymptotically conical, in the following sense.
For every $p\in X$ and every $\eps>0$, there is an $r_1\in[r_0,\infty)$ such that
if $x\in S\setminus B(p,r_1)$, then the angle at $x$ between the geodesic
segment $\ol{xp}$ and $S$ is $<\eps$, and the map 
$S\setminus B(p,r_1)\ra [r_0,\infty)$
given by the distance function from $p$ is a fibration with circle fiber.

3. If the area growth of $S$ is Euclidean, i.e.
$$
\frac{\area(B(p,r)\cap S)}{\pi r^2}\ra 1
$$
as $r\ra \infty$, then $S$ is a $2$-flat.
\end{theorem}
\proof
Using a standard argument, we may assume without loss of generality (and at the cost
of some deterioration in quasi-isometry constants which will be suppressed), 
that $Q$ is the image of a $C$-Lipschitz  $(L,A)$-quasi-isometric embedding $f:\R^2\ra X$,
where $C=C(L,A)$.
The mapping $f$ is proper, and hence induces a homomorphism 
$f_*:H_2^{\lf}(\R^2)\ra H_2^{\lf}(X)$
of locally finite  homology groups.  
We define  $S$ to be the support set of the image of the fundamental class of 
$\R^2$ under $f_*$: 
\begin{equation}
\label{eqnsuppcontainedinimage}
S\defeq\supp\,(f_*([\R^2]))\subset \im(f)=Q.
\end{equation}

\begin{lemma}
There are constants $D=D(L,A)$ and $a=a(L,A)$ such that 
such that:

1. The Hausdorff distance between $S$ and $Q$ is at most $D$.

2. For every 
 $p\in X$, the area of $B(p,r)\cap S$ is at most $a\,(1+r)^2$.
\end{lemma}
\proof
Using the uniform contractibility of $\R^2$, one may construct a proper map 
$g:Q\ra \R^2$ such that $d(g\circ f,\id_{\R^2})$ is bounded by a function
of $(L,A)$.   In particular, the composition of proper maps
$$
\R^2\stackrel{f}{\lra}Q\stackrel{g}{\lra}\R^2
$$
is properly homotopic to $\id_{\R^2}$.  Hence 
 $(g\circ f)_*([\R^2])=[\R^2]$, so $\supp((g\circ f)_*([\R^2]))=\R^2$.
On the other hand 
$$
\supp((g\circ f)_*([\R^2]))\subset \;g(S),
$$
which implies that $Q=\im(f)$
is contained in a controlled neighborhood of $S$.

The last assertion follows from the fact that $S\subset Q$ and
$Q$ has quadratic area growth, being the image of a Lipschitz quasi-isometric embedding.
\qed

Therefore
Theorem \ref{thmquadraticgrowthsupportsetstructure} applies to $S$, and 
by part 2, we get a  fibration 
$$
d_p:S\setminus B(p,r_0)\ra [r_0,\infty)
$$
whose fiber is homeomorphic to a finite graph $\G$  all of whose vertices
have valence $\geq 2$. If $\G$ had a branch point, we could apply
Lemma \ref{lembranching},  contradicting the fact that $S$ is a quasi-flat.
Thus $S$ is locally isometric to $\R^2$ outside $B(p,r_0)$.
\qed

\section{Square complexes}

In this section  $X$ will be a locally finite $\cat(0)$ square complex
with isolated suspensions.   

\bigskip
\begin{remark}
\label{remautomaticisolatedsuspensions}
 It is not difficult to show that if ${\mathcal F}$ is the collection of
$\cat(1)$ graphs $\Ga$ all of whose edges have length $\frac{\pi}{2}$,
then ${\mathcal F}$ has isolated suspensions.  In particular,   any $\cat(0)$
square complex has isolated suspensions.  However, we will not
need this fact for our primary applications, so we omit the proof.
\end{remark}

\begin{theorem}
\label{thmquasiflathalfplanes}
Let $\si\in H_2^{\lf}(X)$ be a quadratic growth locally finite
homology class whose support set $S$ is a quasiflat.  Then there is a
finite collection $\{H_1,\ldots,H_k\}$ of half-plane subcomplexes contained
in $S$, and a finite subcomplex  $W\subset S$
such that 
$$
S\,=\,W\cup\left( \cup_i \;H_i\right)\,.
$$
\end{theorem}
\proof
Pick $p\in X$ and $\eps\in (0,\frac{\pi}{2})$.
Let $r_1$ be as in Theorem \ref{thmqflataflat}, and  set 
$$
Y_1\defeq S\setminus B(p,r_1).
$$ 
Then $Y_1$ is a complete flat Riemannian surface with concave boundary
$\partial Y_1=S(p,r_0)\cap Y_1$.

Now pick $\al\in(0,\frac{\pi}{8})$, $r_2\in [r_1,\infty)$,
and let $Y_2\defeq S\setminus B(p,r_2)$.

\bigskip
\begin{lemma}
\label{lemhalfplane}
Provided $r_2$ is sufficiently large (depending on $\al$),
for every $x\in Y_2$, and every semi-circle $\tau\subset \Si_x\,S$
such that
$$
d(\tau,\log_xp)>\al,
$$
there is a subset $Z\subset S$ isometric to a Euclidean half-plane,
such that $\Si_x\,Z=\tau$.
\end{lemma} 
\proof
First suppose $y\in Y_2$, and $v\in \Si_yS$
is a tangent vector  such that 
$\angle_y(v,\log_yp)>\al$.  Provided $r_2\sin\al>r_1$, 
there will be a unique geodesic ray $\ga_v\subset S$
starting at $y$ with direction $v$; this follows from a continuity
argument, since triangle comparison implies that any geodesic segment
with initial direction $v$ remains outside  $B(p,r_1)$.

If $\tau\subset\Si_xS$ is a semi-circle (i.e. a geodesic segment of length $\pi$), 
and $\angle_x(\tau,\log_xp)>\al$, then the union of the 
rays $\ga_v$, for $v\in \tau$, will form a subset of $S$
isometric to a Euclidean half-plane.  
\qed

\bigskip
\bigskip
{\em Proof of Theorem \ref{thmquasiflathalfplanes} continued.}
We now assume that $r_2$ is large enough that Lemma \ref{lemhalfplane}
applies.

Our next step is to construct a finite collection of half-planes
in $S$.

Consider the boundary $\partial Y_2$.  This is the frontier of the 
set $K\defeq S\cap \ol{B(p,r_2)}$ in $S$.   Since $K$ is locally convex
near $\partial K=\partial Y_2$, it follows that for each $x\in \partial Y_2$, there
is a well-defined space of directions $\Si_xK$, which consists of the 
directions $v\in \Si_xS$ such that $\angle_x(v,\log_xp)\leq \frac{\pi}{2}$.
Also, there is  a normal space $\nu_x \,K\subset \Si_xS$
consisting of the directions $v\in \Si_xS$ making an angle at least
$\frac{\pi}{2}$ with $\Si_xK$.  When $\eps$ is small, the angle
$\angle_x(\log_xp,\Si_xS)$ is small, and hence
$\pi-\angle_x(v,\log_xp)$ will be small for every $v\in \nu_x\,K$.  
In particular, when $\eps$ is small, for every $v\in \nu_vK$
there will be a semi-circle
$\tau_v\subset \Si_xS$ such that 

1. $\tau_v$ makes an angle at least $\frac{\pi}{8}$
with $\log_xp$.

2. If $Z_v\subset S$ is the subset
obtained by applying Lemma
\ref{lemhalfplane} to $\tau_v$, then the boundary of $Z_v$
 is parallel to one of the  sides of a square
$P\subset S$ which contains $x$. 

3. The angle between $\partial Z_v$ and $v$ is at least $\frac{\pi}{8}$.

\noindent
We let $H_v\subset Z_v$ be the largest half-plane subcomplex of $Z_v$.  It
follows from property 2  that
$H_v$ may be obtained from $Z_v$ by removing a strip of thickness
$<1$ around $\partial Z_v$.   

Now let $\H$ be the collection of all  half-planes obtained this
way, where $x$ ranges over $\partial Y_2$, and $v\in \nu_x\,K$.  Observe
that this is a finite collection, since each $H\in \H$  has a boundary
square lying in $B(p,1+r_2)$, and two half-planes $H,H'\in \H$ 
sharing a boundary square must be the same.  

We now claim that 
$$
S\setminus \cup_{H\in\H} \;H
$$
is contained in $\ol{B(p,r_2+\sec\frac{\pi}{8})}$.  To see this note that if 
$y\in Y_2$, then there is a shortest  path in $S$ from $y$ to $K$.  Since $S$ 
is locally convex, this path will 
be a geodesic segment $\ol{yx}$ in $X$, where $x\in \partial Y_2$.  
Let $v\defeq\log_xy\in \Si_x S$.
Then $\ol{yx}$ is contained in $Z_v$, and in view of condition
3 above, all but an initial segment
of length at most $\sec\frac{\pi}{8}$ will be contained in $H_v\subset Z_v$.
The claim follows.  
\qed

\bibliography{quasiflat}

\begin{thebibliography}{BKMM12}

\bibitem[BH99]{bridsonhaefliger}
M.R. Bridson and A.~Haefliger.
\newblock {\em Metric spaces of non-positive curvature}.
\newblock Springer-Verlag, Berlin, 1999.

\bibitem[BKMM12]{behrstockkleinerminskymosher}
J.~Behrstock, B.~Kleiner, Y.~Minsky, and L.~Mosher.
\newblock Geometry and rigidity of mapping class groups.
\newblock {\em Geom. Topol.}, 16(2):781--888, 2012.

\bibitem[BKS08]{raag}
M.~Bestvina, B.~Kleiner, and M.~Sageev.
\newblock The asymptotic geometry of right-angled {A}rtin groups. {I}.
\newblock {\em Geom. Topol.}, 12(3):1653--1699, 2008.

\bibitem[EF97]{eskinfarb}
A.~Eskin and B.~Farb.
\newblock Quasi-flats and rigidity in higher rank symmetric spaces.
\newblock {\em J. Amer. Math. Soc.}, 10(3):653--692, 1997.

\bibitem[Esk98]{eskin}
A.~Eskin.
\newblock Quasi-isometric rigidity of nonuniform lattices in higher rank
  symmetric spaces.
\newblock {\em J. Amer. Math. Soc.}, 11(2):321--361, 1998.

\bibitem[Hat02]{hatcher}
Allen Hatcher.
\newblock {\em Algebraic topology}.
\newblock Cambridge University Press, Cambridge, 2002.

\bibitem[KL]{klla}
B.~Kleiner and U.~Lang.
\newblock Quasi-minimizers in hadamard spaces.
\newblock in preparation.

\bibitem[KL97a]{kaple}
M.~Kapovich and B.~Leeb.
\newblock Quasi-isometries preserve the geometric decomposition of {H}aken
  manifolds.
\newblock {\em Invent. Math.}, 128(2):393--416, 1997.

\bibitem[KL97b]{klle}
B.~Kleiner and B.~Leeb.
\newblock Rigidity of quasi-isometries for symmetric spaces and {E}uclidean
  buildings.
\newblock {\em Inst. Hautes \'Etudes Sci. Publ. Math.}, (86):115--197, 1997.

\bibitem[Mos73]{mostow}
G.~D. Mostow.
\newblock {\em Strong rigidity of locally symmetric spaces}.
\newblock Princeton University Press, Princeton, N.J., 1973.
\newblock Annals of Mathematics Studies, No. 78.

\bibitem[Xie05]{xie}
X.~Xie.
\newblock The {T}its boundary of a {$\rm CAT(0)$} 2-complex.
\newblock {\em Trans. Amer. Math. Soc.}, 357(4):1627--1661 (electronic), 2005.

\end{thebibliography}
\bibliographystyle{alpha}

\end{document}